# COMPACTLY SUPPORTED WAVELETS DERIVED FROM LEGENDRE POLYNOMIALS: SPHERICAL HARMONIC WAVELETS

M.M.S. LIRA[1], H.M. DE OLIVEIRA[2], M.A. CARVALHO JR[1], R.M.C. DE SOUZA[2]

Federal University of Pernambuco - UFPE
LDSP - Power Systems Digital Laboratory[1]; CODEC - Communications Research Group[2]
C.P. 7800, ZIP 50711-970, Recife - PE, Brazil
E-mail: {milde, hmo, macj, ricardo}@ufpe.br   http://www.ee.ufpe.br/codec/Legendre.html

*Abstract:* A new family of wavelets is introduced, which is associated with Legendre polynomials. These wavelets, termed spherical harmonic or Legendre wavelets, possess compact support. The method for the wavelet construction is derived from the association of ordinary second order differential equations with multiresolution filters. The low-pass filter associated to Legendre multiresolution analysis is a linear phase finite impulse response filter (FIR).

## 1 Introduction

A recent paper [1] introduced new orthogonal wavelets related to the Mathieu wave equation [2]. It was then suggested the idea of linking other ordinary 2$^{nd}$-order differential equations (wave equations) with the transfer function of analysing orthogonal multiresolution filters [3]. An appealing case comprises orthogonal polynomials that are solution of certain differential equations. Different systems of orthogonal polynomials satisfying second-order ordinary differential equations have long been used in Physics and Engineering [4]. Jacobi, Legendre, Hermite, Gegenbauer, Laguerre are examples of polynomials encompassing very well-known properties [2, 5].

 This paper is particularly concerned with Legendre polynomials. Associated Legendre polynomials are the colatitudinal part of the spherical harmonics which are common to all separations of Laplace's equation in spherical polar coordinates. The radial part of the solution varies from one potential to another, but the harmonics are always the same and are a consequence of spherical symmetry. Legendre functions have widespread applications including Control Systems, Engineering and Physics [4]. They are used in problems in which spherical polar are appropriate, for instance, to determine the steady-state temperature in a uniform sphere of radius unity when one half of the surface is kept at 0$^{o}$C and the other half at 1$^{o}$C [4]. Further attractive matter defined on a spherical surface can be found in electromagnetics and acoustics [6, 7]. They have also been useful in system identification [8], as a method for solving the least-square estimation problem [9], and even in sophisticated imaging problems such as "shape from shading algorithms" for 3-D image reconstruction [10]. Legendre polynomials have recently been used so as to define a windowing technique for FIR filter design, instead of Hamming, Kaiser or other standard window [11]. They are also useful for analyzing the propagation characteristic of multimode planar graded-index optical fibers [12] or for investigating scattering and waveguiding by planar periodic strip [13].

 Transform coding is now a very common tool in signal and image analysis. Discrete versions (DWT) of wavelet transforms have successfully been proposed to be used instead of classical discrete transforms such as Discrete Fourier Transform (DFT) and Discrete Cosine Transform (DCT), in several applications [14-16]. Another discrete transform less recognized is the discrete Legendre transform (DLT) [17]. They are especially efficient when data can be modeled by a polynomial function [7]. Parallel algorithms for computing DLT and IDLT are available [17,18]. Colomer & Colomer recently introduced a new electrocardiogram (ECG) compression scheme, which is based on the DLT [19]. The DLT can also be a powerful tool in computerized tomography [20]. Wavelets associated to finite impulse response filters (FIR) are commonly preferred in most applications.

Each signal (waveform signal) possesses infinitely many different representations. There exists as well an infinite number of wavelet systems whose applicability and relevance strongly depends on the analysing signal properties. Besides conventional wavelets such as Morlet, Meyer, Battle, Daubechies, Coiflets etc. [21], other wavelets have recently been proposed in the literature [e.g. 1, 22-23]. It seems to be natural to look for new wavelets holding further particular symmetries. We are led to ask: Are there simple wavelets associated to linear phase FIR? Are there Legendre wavelets? If the answer is affirmative, then how DWT based on Legendre functions can be compared with the DLT? This note answers the first question and furnishes the basis to investigate the second one.

## 2 Preliminaries On Legendre Functions

Spherical harmonics are solutions of the Legendre 2nd-order differential equation, $n$ integer:

$$(1-z^2)\frac{d^2y}{dz^2} - 2z\frac{dy}{dz} + n(n+1)y = 0. \qquad (1)$$

The solution of (1) above is the $n^{th}$ order Legendre polynomial $P_n(z)$. Such polynomials can easily be found from the Rodrigues' formula [2]:

$$P_n(z) = \frac{1}{2^n n!}\frac{d^n}{dz^n}\left(z^2-1\right)^n. \qquad (2)$$

Frequently it is more appropriate to confine the variable $z$ to be real, $z=x$ into the region $|x|\leq 1$. In such cases, it is common to adopt the variable change $x=\cos\theta$ and deal with the polynomials under the form $P_n(\cos\theta)$. It follows that $P_n$ is $2\pi$-periodic and $|P_n(\cos\theta)|\leq 1$. This is precisely what we require, assuming $\theta$ related to the spectral frequency $\omega$. Additionally, the following property holds:

$$P_n(-x) = (-1)^n P_n(x)$$
$$\text{so that } P_{2n+1}(-x) = -P_{2n+1}(x). \qquad (3)$$

The $n^{th}$ Legendre polynomial has $n$ distinct roots within the interval $-1\leq x\leq 1$.

$|P_n(\cos\theta)|$ functions present a shape that visibly resemble that of a low-pass filter, provided that $n$ is odd [2]. Limiting values are $P_n(1)=1$ and $P_n(0)=0$, $n$ odd. The Legendre polynomials are described by:

$$P_1(\cos\theta) = \cos\theta$$

$$P_3(\cos\theta) = \frac{1}{8}(5\cos3\theta + 3\cos\theta)$$

$$P_5(\cos\theta) = \frac{1}{128}(63\cos5\theta + 35\cos3\theta + 30\cos\theta)$$

$$\ldots \qquad (4)$$

## 3 Legendre Multiresolution Filters

$P_n(\cos\theta)$ polynomials can be used to define the smoothing filter $H(\omega)$ of a multiresolution analysis (MRA) [3]. Since the appropriate boundary conditions for an MRA are $|H(0)|=1$ and $|H(\pi)|=0$, the smoothing filter of an MRA can be defined so that the magnitude of the low-pass $|H(\omega)|$ can be associated to Legendre polynomials according to:

$$|H(\omega)| = |P_{2n+1}\left(\cos\frac{\omega}{2}\right)|. \qquad (5)$$

This is quite similar to the approach implemented to derive the Mathieu wavelet [1]: The cosine-elliptic Mathieu function was associated to a low-pass filter according to

$$|H_\nu(\omega)| = \left|\frac{ce_\nu\left(\frac{\omega}{2},q\right)}{ce_\nu(0,q)}\right|. \qquad (6)$$

Here, $\nu$ plays the role of $2n+1$: Mathieu wavelets can be only derived when $\nu$ is odd. Moreover Legendre polynomials are already normalized, i.e., $|P_{2n+1}(\cos 0)| = |P_{2n+1}(1)| = 1$. Henceforth, we assume $\nu=2n+1$. The magnitude of the low-pass is given by

$$|H_\nu(\omega)| = \left|\frac{P_\nu\left(\cos\frac{\omega}{2}\right)}{P_\nu(\cos 0)}\right|. \qquad (7)$$

Illustrative examples of filter transfer functions for a Legendre MRA are shown in figure 1, for $\nu=1,3$ and 5. We can notice the low pass behaviour for the filter $H$, as expected. The number of zeroes within $-\pi<w\leq\pi$ is equal to the degree of the Legendre polynomial. Therefore, the roll-off of side-lobes with frequency is easily controlled by the parameter $\nu$.

A suitable phase assignment should be done so as to properly adjust the transfer function $H_\nu(w)$ to the form

$$H_\nu(w) = \frac{1}{\sqrt{2}}\sum_{k\in Z}h_k^\nu e^{-jwk}. \qquad (8)$$

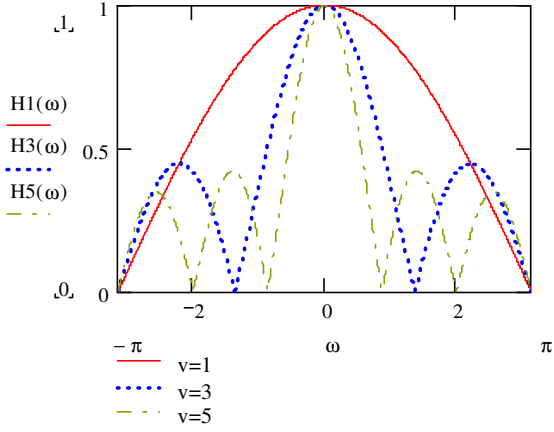

Fig. 1 - Magnitude of the transfer function for Legendre multiresolution smoothing filters. Filter $|H_\nu(\omega)|$ for a few orders: $\nu=1$ (solid line), $\nu=3$ (dot line), and $\nu=5$ (dashdot line).

In Mathieu's case, the assumption was [1]

$$H_\nu(\omega) = -e^{-j\nu\frac{\omega}{2}} \frac{ce_\nu\left(\frac{\omega}{2}, q\right)}{ce_\nu(0, q)}. \quad (9)$$

A possible proper and rather similar solution takes for granted that

$$H_\nu(\omega) = -e^{-j\nu\frac{\omega}{2}} P_\nu\left(\cos\frac{\omega}{2}\right). \quad (10)$$

An auxiliary function $y_\nu(\omega)$ defined by

$$y_\nu(\omega) = -e^{j\nu\omega} H_\nu(2\omega) \quad (11)$$

satisfies a Legendre differential equation.

In order to derive a multiresolution analysis, the transfer function of the high pass analysing filter $G_\nu(\omega)$ can be chosen applying the QMF condition [3, 21], yielding:

$$G_\nu(\omega) = e^{j(\nu-2)\frac{(\omega-\pi)}{2}} P_\nu\left(\sin\frac{\omega}{2}\right). \quad (12)$$

It can immediately be seen that $|G_\nu(0)|=0$ and $|G_\nu(\pi)|=1$, as expected.

Once defined $H_\nu(\omega)$, the next step is to compute the filter coefficients $\{h_k\}$, $k \in Z$. This can be done by applying explicit expressions involving Legendre polynomials and trigonometric functions [2,5]:

$$P_n(\cos\theta) = \sum_{m=0}^{n} a_m \cos((n-2m)\theta), \quad 0 < \theta \leq \pi, \quad (13)$$

where $a_m = \frac{1}{4^n}\binom{2m}{m}\binom{2n-2m}{n-m}$.

Therefore, substituting (13) into (10)

$$H_\nu(\omega) = -e^{-j\nu\frac{\omega}{2}} \sum_{m=0}^{\nu} a_m \cos\left[(\nu-2m)\frac{\omega}{2}\right]. \quad (14)$$

After a rather trivia manipulation

$$H_\nu(\omega) = \sum_{m=0}^{\nu} -\frac{a_m + a_{\nu-m}}{2} e^{-jm\omega}. \quad (15)$$

A term by term comparison between (8) and (15) implies

$$\frac{h_k^\nu}{\sqrt{2}} = -\frac{a_k + a_{\nu-k}}{2}, \quad k=0,1,2,...,\nu. \quad (16)$$

It follows then the symmetry: $h_k^\nu = h_{\nu-k}^\nu$. There are just $\nu+1$ non-zero filter coefficients on $H_\nu(\omega)$, so that the Legendre wavelets have compact support for every odd integer $\nu$. Finally,

$$\frac{h_k^\nu}{\sqrt{2}} = -\frac{1}{2^{2\nu}}\left[\binom{2k}{k}\binom{2\nu-2k}{\nu-k}\right]. \quad (17)$$

Table 1. Smoothing Legendre FIR filter coefficients for $\nu=1,3,5$ ($N$ is the wavelet order.)

|     | $\nu=1$ ($N=1$) | $\nu=3$ ($N=2$) | $\nu=5$ ($N=3$) |
|-----|-----|-----|-----|
| $h_0$ | $-\frac{\sqrt{2}}{2}$ | $-\frac{5\sqrt{2}}{16}$ | $-\frac{63\sqrt{2}}{256}$ |
| $h_1$ | $-\frac{\sqrt{2}}{2}$ | $-\frac{3\sqrt{2}}{16}$ | $-\frac{35\sqrt{2}}{256}$ |
| $h_2$ |  | $-\frac{3\sqrt{2}}{16}$ | $-\frac{30\sqrt{2}}{256}$ |
| $h_3$ |  | $-\frac{5\sqrt{2}}{16}$ | $-\frac{30\sqrt{2}}{256}$ |
| $h_4$ |  |  | $-\frac{35\sqrt{2}}{256}$ |
| $h_5$ |  |  | $-\frac{63\sqrt{2}}{256}$ |

N.B. The minus signal can be suppressed.

Legendre wavelets can be derived from the low-pass reconstruction filter by an iterative procedure (the cascade algorithm). The wavelet has compact support and finite impulse response AMR filters (FIR) are used (table 1). Curiously, the first wavelet of the Legendre's family is exactly the well-known Haar wavelet. Figures 2 and 3 show an emerging pattern that progressively looks like the wavelet's shape. As with many wavelets there is no nice analytical formula for describing harmonic spherical wavelets.

Aiming to investigate potential applications of such wavelets, software to compute them should be written. Nowadays one of the most powerful software supporting wavelet analysis is the Matlab[TM] (MATrix LABoratory) [24], especially when the wavelet graphic interface is available. In the wavelet toolbox, there exists five kinds of wavelets (type the command *waveinfo* on the prompt): (*i*) crude wavelets (*ii*) Infinitely regular wavelets (*iii*) Orthogonal and compactly supported wavelets (*iv*) biorthogonal and compactly supported wavelet pairs. (*v*) complex wavelets. The majority of functional wavelets are of kind three or four.

Legendre wavelets can be easily loaded into the MatLab™ wavelet toolbox -- The *m*-files to allow the computation of Legendre wavelet transform, details and filter are currently (*freeware*) available at the URL: http://www.ee.ufpe.br/codec/Legendre.html
The Legendre wavelet shape can easily be visualised using the *wavenemu* command of Matlab™.

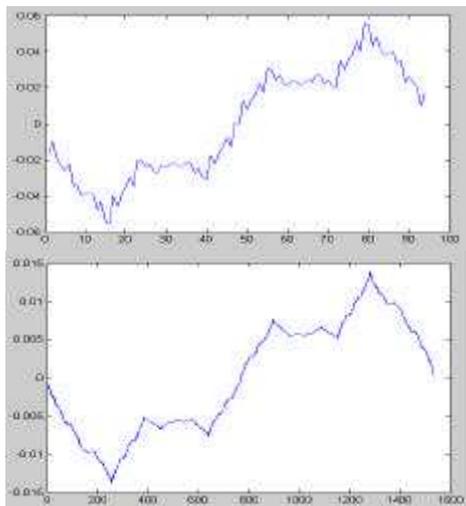

Fig. 2 - Shape of Legendre Wavelets of degree ν=3 (legd2) derived after 4 and 8 iteration of the cascade algorithm, respectively.

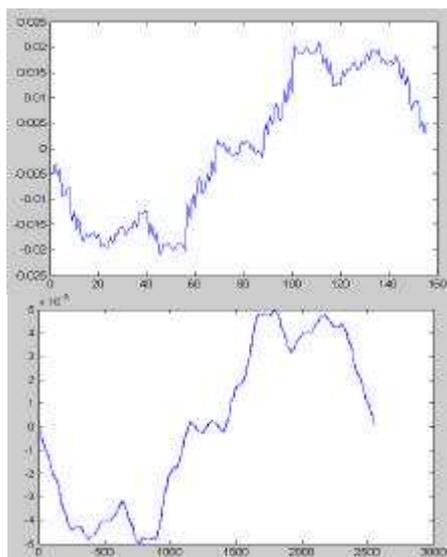

Fig. 3 - Shape of Legendre Wavelets of degree ν=5 (legd3) derived by the cascade algorithm after 4 and 8 iterations of the cascade algorithm, respectively.

The finite support width Legendre family is denoted by legd (short name). Wavelets: 'legd$N$'. The parameter $N$ in the legd$N$ family is found according to $2N=\nu+1$ (length of the MRA filters).

Figure 4 shows legd8 wavelet display using MatLab™. The association between Legendre and other linear phase wavelets [26] remained to be investigated.

As with many other wavelets, the filter coefficients cannot guarantee a perfect reconstruction [27]. However, starting with an non-orthogonal compactly supported wavelet basis, it is possible to generate an orthogonal multiresolution by means of Herley-Vetterli approach [28]. Furthermore, Biorthogonal Legendre wavelets could probably be derived to achieve exact reconstruction property. This seem to be the most promising approach for further development.

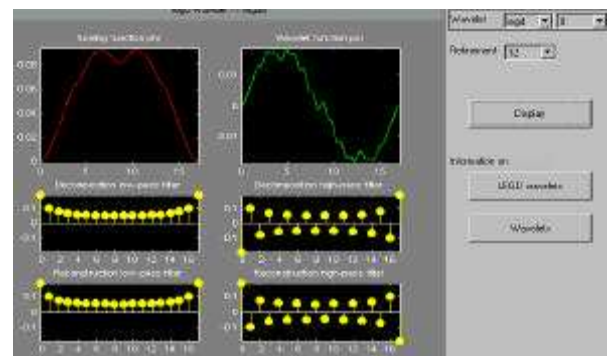

Fig. 4 - legd8 wavelet display over Matlab™ using the *wavemenu* command.

2-dimensional Legendre transform (2D-legd$N$) for image analysis is promptly available. Interesting effects have been observed in multiresolution decomposition of hand-drawn images such as illustrated in figure 5.

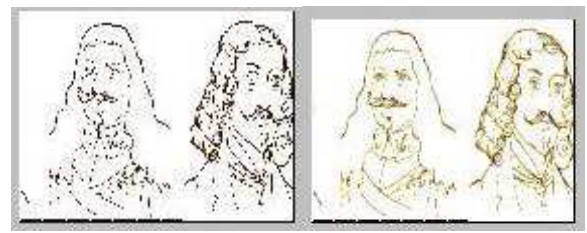

Fig. 5 - Reconstruction of a hand-drawn draft using 2D-legd2. (a) Original (b) Reconstructed image from 1-level wavelet decomposition.

Wavelet Packets (WP) systems derived from Legendre wavelets can also be easily accomplished. Figure 6 illustrates the WP functions derived from legd2.

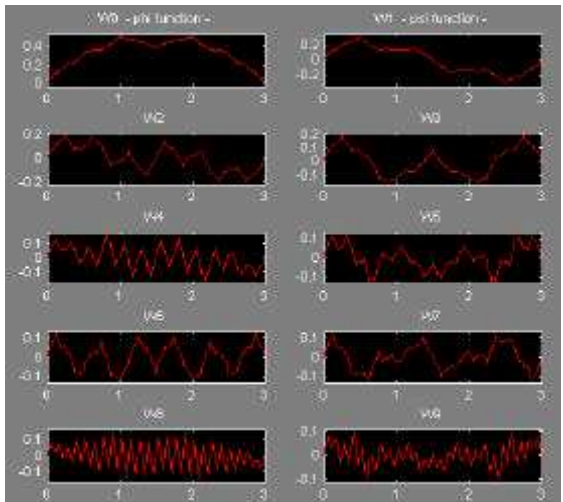

Fig. 6 - Legendre (legd2) Wavelet Packets W system functions: WP from 0 to 9.

## 4 Final Remarks

A family of harmonic spherical wavelets was introduced. It is shown that the transfer functions of the corresponding multiresolution filters are related to Legendre polynomials. In contrast with Mathieu wavelets derived from a similar method [1], the Legendre wavelets have compact support, which is very attractive from the practical viewpoint. The magnitude of the smoothing filters corresponds to Legendre functions with odd order. An extra appealing feature is that the Legendre filters are *linear phase* FIR (i.e. multiresolution analysis associated to linear-phase filters). Although being compactly supported wavelet, legd*N* are not orthogonal (except for *N*=1). Dual scaling and dual filter are currently being designed to build biorthogonal wavelets. Any non-orthogonal compactly supported wavelet analysis can be replaced by an orthogonal one, which is implemented by (realisable) IIR filter. Further generalisations should be examined, e.g. using $P_n^m(z)$ associate Legendre functions or using even order functions. These wavelets have been implemented on Matlab$^{TM}$ (wavelet toolbox). The relevance of such new wavelets on signal analysis is currently being investigated in our laboratory. Potential application of such wavelets includes imaging, electromagnetics, optics, acoustics, and electrocardiogram (ECG) data compression, among other topics.

*Acknowledgements:*
The authors express their indebtedness to Mr. R.J.S. Cintra and Mrs L.R. Soares (UFPE, Brazil) for helpful comments and assistance.